\documentclass{amsart}
\usepackage{amsmath}
\usepackage{amssymb}
\usepackage{graphicx}
\usepackage{latexsym}
\usepackage{color}
\usepackage{amscd}
\usepackage[all]{xy}
\usepackage{enumerate}
\usepackage{hyperref}
\usepackage{subfigure}
\usepackage{soul}
\newtheorem{thm}{Theorem}

\newtheorem{theorem}[thm]{Theorem}

\theoremstyle{definition}

\newtheorem*{definition*}{Definition}

\newtheorem{remark}[thm]{Remark}

\newcommand{\CPb}{\overline{\mathbb{CP}}{}^{2}}
\newcommand{\CP}{{\mathbb{CP}}{}^{2}}

\newcommand{\Z}{\mathbb{Z}}

\begin{document}

\title[Dissolving knot surgered $4$--manifolds by classical cobordism arguments]
{Dissolving knot surgered $4$--manifolds \\ by classical cobordism arguments}

\author[R. \.{I}. Baykur]{R. \.{I}nan\c{c} Baykur}
\address{Department of Mathematics and Statistics, University of Massachusetts, Amherst, MA 01003-9305, USA}
\email{baykur@math.umass.edu}

\begin{abstract}
The purpose of this note is to show that classical cobordism arguments, which go back to the pioneering works of Mandelbaum and Moishezon, provide quick and unified proofs of any knot surgered compact simply-connected $4$--manifold $X_K$ becoming diffeomorphic to $X$ after a single stabilization by connected summing with $S^2 \times S^2$ or $\CP \# \CPb$, and almost complete decomposability of $X_K$ for many almost completely decomposable $X$, such as the elliptic surfaces.
\end{abstract}

\maketitle

\setcounter{secnumdepth}{2}
\setcounter{section}{0}

In the late 1970s, Moishezon and Mandelbaum discovered that many algebraic surfaces, such as the elliptic surfaces $E(n)$, are almost completely decomposable, i.e. after a single connected sum with the complex projective plane $\CP$, they \emph{dissolve} --become diffeomorphic-- to a connected sum of $\CP$s and $\CPb$s; see for example \cite{Man0, Man1, Man2, ManMos1, ManMos2, Mos}. Most recently, Choi, Park and Yun considered \emph{knot surgered elliptic surfaces} $E(n)_K$ of Fintushel and Stern \cite{FS}, which (for each $n$) constitute a large family of irreducible $4$--manifolds that contain infinitely many symplectic but not algebraic, and even non-symplectic  \mbox{$4$--manifolds} that are homeomorphic but not pairwise diffeomorphic \cite{FS}. Using handle diagrams and Kirby calculus, they proved that, after connected summing with a single $\CP$, all $E(n)_K$ become diffeomorphic to each other (in particular to $E(n)$ when $K$ is the unknot), and  concluded that $E(n)_K$ are ACD  \cite{CPY}. Here we give an alternate proof  using $5$--dimensional cobordism arguments, which follow the classical framework of Mandelbaum and Moishezon ---with further contributions from Gompf \cite{Gompf0}. These general arguments show that many other knot surgered $4$--manifolds  are almost completely decomposable (e.g. many irreducible $4$--manifolds in \cite{JPark}). These are knot surgery \mbox{$4$--manifolds} that can be seen as fiber sums of almost completely decomposable $4$-manifolds and knot surgered $E(n)$).) A simpler version of our proof provides general results on the stable equivalence of knot surgered $4$--manifolds after a single stabilization by connected summing with $S^2 \times S^2$ or $\CP \# \CPb$.

\smallskip
\begin{theorem} \label{mainthm}
Let $X$ be a compact, simply-connected, smooth $4$--manifold, and  $T$ be a square zero torus in $X$, with  simply-connected complement. For a knot $K$ in $S^3$, let $X_K$ be the simply-connected $4$--manifold obtained by knot surgery along $T$ in $X$. Then, $X_K \# \,  S^2 \times S^2 \cong X \# \, S^2 \times S^2$, and when $X$ is not spin,  \mbox{$X_K \, \# \CP \# \CPb  \cong X \, \# \CP \# \CPb$} as well. Moreover, if $(X,T)$ is obtained from $(X',T')$ by blowing-up at a point on $T$, then we have $X_K \# \, \CP \cong X' \# \, \CP \# \CPb$, and therefore $X_K$ is almost completely decomposable whenever $X'$ is. It also follows that any knot surgered elliptic surface $E(n)_K$ is almost completely decomposable.
\end{theorem}

\vspace{-0.15cm}
The first part of our theorem recaptures somewhat more specific results of Auckly and Akbulut, who used handle diagrams and Kirby calculus to show that knot surgered simply-connected $4$--manifolds $X_K$ become diffeomorphic to $X$ after a stabilization by connected summing with $\CP \# \CPb$ or with $S^2 \times S^2$, provided the torus fiber is contained in a Gompf nucleus \cite{Auckly}, or in a cusp neighborhood \cite{Akbulut}, respectively. Other general results on the sufficiency of single stabilizations for homeomorphic simply-connected $4$--manifolds to become diffeomorphic for various constructions of exotic $4$--manifolds were given by the author and Sunukjian in \cite{BS} using \mbox{$5$--dimensional} \emph{round handle cobordisms.}  

\vspace{0.15in}
Below $(X_1, F_1) \#_\phi (X_2, F_2):= (X_1 \setminus \nu F_1) \cup_{\phi} (X_2 \setminus \nu F_2$) denotes the \emph{(generalized) fiber sum} \cite{Gompf} of $X_i$, $i=1,2$, along square zero closed orientable surfaces $F_i \subset X_i$ of the same genera, with the gluing map $\phi$. Let $K$ be a knot in $S^3$. We can recast the construction of $X_K$ as a fiber sum: for $T$ a square zero torus with simply-connected complement in $X$, we have $X_K = (X, T) \#_{\phi} (S^1 \times S^3, T_K)$, where \mbox{$T_K = S^1 \times K$.} The gluing map $\phi$ identifies the two circles $a$ and $b$, which frame the $T$ factor of $T \times D^2 \cong \nu T$, with $S^1 \times \text{pt}$ and $\text{pt} \times K$, which frame the $T_K$ factor of $T_K \times D^2 \cong \nu \, T_K= S^1 \times \nu K$. (Here $\phi$ identifies the meridian of $T$, which is null-homotopic in $X \setminus \, \nu T$ with $\text{pt} \times \mu$, where $\mu$ is the meridian of $K$, normally generating $\pi_1(S^3 \setminus \nu K)$.) When $X$ is the elliptic surface $E(n)$, $T$ is a regular elliptic fiber \cite{FS}.

\begin{proof} 
There is a compact, oriented cobordism $W$ from the disjoint union $X_1 \sqcup X_2$ to $Z=(X_1, F_1) \#_\phi (X_2, F_2)$ obtained by gluing a $D^1 \times \Sigma_g \times D^2$ (for $F_i \cong \Sigma_g)$ to \mbox{$D^1 \times (X_1 \sqcup X_2)$,} where $\partial D^1 \times \Sigma_g \times D^2$ is identified with $\nu F_1 \,\sqcup \, \nu F_2$ using $id \,\sqcup \, \phi$. We can break down this cobordism by gluing $D^1 \times \Sigma_g \times D^2$ in pieces, using the standard handle decomposition on $\Sigma_g=h^0+\sum_{i=1}^{2g} h^1_i \ + h^2$. After gluing $D^1 \times h^0 \times D^2$, we obtain $X_1 \# X_2$, the connected sum of the two manifolds. Observe that the gluing of $D^1 \times h^1_i \times D^2$ is the same as attaching a $5$--dimensional $2$-handle $H_i$. Looking at the cobordism upside down, we see that the gluing of $D^1 \times h^2 \times D^2$ is the same as attaching a $5$--dimensional $2$--handle as well, but to $Z$ instead. If a $5$--dimensional $2$--handle is attached to a simply-connected non-spin $4$--manifold, let it be $Z$, $X_1$ or $X_2$, the result is a connected sum with $\CP \# \CPb$, or if we like, with $S^2 \times S^2$, since the framings can be chosen either way in this case. If all are simply-connected, the second to last level of the cobordism shows that $X_1 \# X_2 \# 2g (\CP \# \CPb) \cong Z \# \CP \# \CPb$.  (On the other hand, if $Z, X_1, X_2$ are spin, the framings are chosen so that  $X_1 \# X_2 \# 2g  \ S^2 \times S^2 \cong Z \# \, S^2 \times S^2$; perhaps most explicitly stated in \cite{Gompf0}[Lemma 4].) Now, the key observation of Mandelbaum is the following: If $(X_1, F_1)$ is obtained from $(X'_1, F'_1)$ after a blow-up on a point on $F'_1$, then one can trace this exceptional sphere through the cobordism $W$ and eliminate it to arrive at a finer conclusion; namely, $X'_1 \# X_2 \# 2g (\CP \# \CPb) \cong Z \# \CP$ \cite{Man1}[Theorem B\,(2)]. More generally, if $X_2$ is not simply-connected, what we get is a $4$--manifold obtained from $X'_1 \# X_2$ after $2g$ $5$--dimensional $2$--handle attachments along $2g$ disjoint circles in $X_2$, diffeomorphic to $Z \# \CP$ \cite{Man2}[Theorem 2.6\,(ii)])

Now for knot surgered $4$--manifolds, let us first consider the most specific case of \mbox{$X_1=E(1)$} and $Z=E(1)_K$.  We note that the following is essentially a generalization of the classical proof of the almost complete decomposability for $E(n)$, since $E(n)=E(n)_K$, for $K$ the unknot.  Because $E(n)_K \cong (E(1)_K, T) \# (E(n-1), T)$, the arguments in the previous paragraph show that it suffices to prove the almost complete decomposability for $E(1)_K$. 
Observe that $(E(1), F)$ is obtained from $(\CP \# \, 8 \CPb, F')$, where $F'$ is an $8$ times blown-up cubic in $\CP$.  Now for  $E(1)_K = (E(1), T) \#_{\phi} (S^1 \times S^3, T_K)$, the only difference is that \mbox{$S^1 \times S^3$} is not simply-connected. 
In this case, we attach two $5$--dimensional $2$--handles $H_1$ and $H_2$ to \mbox{$\CP \# \,8 \CPb \# \,S^1 \times S^3$} to arrive at $E(1)_K \# \,\CP$. By the chosen gluing map $\phi$, the attaching circle of $H_1$ (resp. $H_2$) is the connected sum of $a$ (resp. $b$) in $T' \subset \CP \# \,8 \CPb$ and $S^1 \times \text{pt}$ (resp. $\text{pt} \times K$) in $T_K \subset S^1 \times S^3$. As $a$ and $b$ are null-homotopic in $\CP \# \,8 \CPb \setminus \nu T'$, the attaching circles of $H_1, H_2$ are isotopic to the latter summands. (Homotopy implies isotopy for $1$--manifolds in a $4$--manifold.) Attaching $H_1$  along $S^1 \times \text{pt}$ of $S^1 \times S^3$ yields $S^4$. (To see this, consider for example the standard handle decomposition of $S^1 \times S^3$ with a single handle for every index but $2$. We then trade the $1$--handle with a $0$--framed $2$--handle attached along an unknot, which cancels against the $3$--handle.) Now, attaching $H_2$ along ---the image of--- $\text{pt} \times K$ in $S^4$ amounts to a connected sum with $S^2 \times S^2$ or $\CP \# \CPb$ (which doesn't matter here, since the other summand $\CP \# \,8 \CPb$ is non-spin). Hence, the second to last level of the cobordism reads from the bottom as $\CP \# 8 \CPb \# S^2 \times S^2 \cong 2\CP \# 9  \CPb$. From the top, it reads as $E(1)_K \# \CP$, showing that $E(1)_K$ is ACD. This in turn concludes that $E(n)_K$ is ACD.   

At this point it should be clear that running the above arguments for arbitrary $X=X_1$ instead of $E(1)$, under respective assumptions on $X$, provides a proof of the second part of the theorem. For the stabilization results in the first part, we run the same arguments simply without eliminating the extra $\CP$, as we discussed in the very first paragraph of our proof.
\end{proof}

\smallskip
\begin{remark}
One can regard the knot surgery operation also  as the fiber sum  $X_K= (X, T) \#_{\phi} (S^1 \times Y_K, T_{\mu})$,  where $Y_K$ is the $3$--manifold obtained by $0$--surgery along the knot $K$ in $S^3$, and $T_\mu=S^1 \times \mu$, for $\mu$ the image of the meridian of $K$. This is a better view point when one would like to perform knot surgery \emph{symplectically}: if $K$ is a genus--$g$ fibered knot, then $S^1 \times Y_K$ admits a symplectic  genus--$g$ fibration over $T^2$, where $T_\mu$ is a symplectic section. Thus, for $T$ a symplectic torus in a symplectic $4$--manifolds $X$, we can take a symplectic fiber sum to arrive at a symplectic \mbox{$X_K$ \cite{FS}.} Note that, although $S^1 \times S^3 \setminus \, \nu T_K  = S^1 \times (S^3 \setminus \nu K)= S^1 \times Y_K \setminus \nu \, T_\mu $, so the knot surgery operation sews the same complements to $X \setminus \nu T$ to produce $X_K$, the roles of $K$ and $\mu$ in $\partial T_K$ and $\partial T_\mu$ are switched. In the latter, the gluing map $\phi$ identifies the two circles $a$ and $b$, which frame the $T$ factor of $T \times D^2 \cong \nu \, T$, with $S^1 \times \text{pt}$ and $\text{pt} \times \mu$, which frame the $T_\mu$ factor of $T_\mu \times D^2 \cong \nu \, T_\mu$. 


With these in mind, we can run similar arguments as in the proof of Theorem~\ref{mainthm}, now for a cobordism $W$ from the disjoint union $X \sqcup \, S^1 \times Y_K$ to $X_K$. In fact, we can make it very comparable to our earlier cobordism if we regard $S^1 \times Y_K$ as the \mbox{$4$--manifold} obtained from $S^1 \times S^3$ by attaching a \mbox{$5$--dimensional} \emph{round $2$--handle} along $S^1 \times K$, which consists of attaching a $5$--dimensional $2$--handle $H^2$ (along $\text{pt} \times K$ in $S^1 \times S^3$) and a $5$--dimensional $3$--handle $H^3$ that goes over $H^2$ geometrically twice but algebraically zero times; see e.g. \cite{BS}. This way, we can extend the cobordism $W$ on one end to once again consider a cobordism $W'$ from $X \sqcup \, S^1 \times S^3)$ to $X_K$. Now, the $3$--handle ($D^1 \times h^2 \times D^2$) of the original cobordism $W$ is attached along a $2$--sphere that is the union of a normal disk to $T$, a normal disk to $\mu$ (which is a $2$--disk capping $K$ in $Y_K$), and a cylinder between the meridian $\mu_T$ of $T$ and $K$. So it goes over the  
$2$--handle $H^2$ of the round $2$--handle once, and we can cancel them. Attaching the $1$--handle $D^1 \times h^0 \times D^2$, and the $2$--handles $H_1$ and $H_2$ to $X \sqcup \, S^1 \times S^3$, we get at a connected sum of $X$ with $S^2 \times S^2$ or $\CP \# \CPb$ as before. Attaching the $3$--handle $H^3$ of the round $2$--handle to the other end $X_K$ as a $2$--handle, we arrive at the same conclusions.
\end{remark}

\begin{remark} Knot surgery using fibered knots, generalized logarithmic transforms (sometimes as Luttinger surgeries along Lagrangian tori), and generalized fiber sums have been the most widely employed techniques for constructing new symplectic $4$--manifolds over the past three decades (see e.g. the  survey \cite{FS2}). Similar cobordism arguments provide analogous results for $4$--manifolds obtained by \emph{(generalized) logarithmic transforms}; see e.g. \cite{FS2} for a definition. For $T$ a square zero torus in $X$ with simply-connected complement,  let $X'$ denote the simply-connected $4$--manifold obtained by a logarithmic transform along $T$ in $X$. We can recast this construction as a fiber sum, too: $X' = (X \setminus \nu T) \cup _{\phi}  (T^2 \times S^2 \setminus  \nu T'')$, where $T''=T^2 \times \text{pt}$. This is the same as the case of  knot surgery with $K$ the unknot, except that the gluing map $\phi$ is chosen to identify the meridian of $T''$ with a primitive $(p,q, r)$ curve in $\Z^3 \cong \partial \nu T$. Although varying these choices for $\phi$ may vary the diffeomorphism type of $X'$, it makes no difference for our cobordism arguments; $5$--dimensional $2$--handle attachments  prescribed by $\phi$ still kill the two generators of $\pi_1(T^2 \times S^2)$ by $0$--surgeries along $S^1 \times \text{pt} \times \text{pt}$ and $\text{pt} \times S^1 \times \text{pt}$. As our discussions illustrate,  almost complete decomposability of most of these $4$--manifolds can be reduced to almost complete decomposability of the simply-connected``building blocks'', which are often algebraic surfaces themselves ---as in the original works of Mandelbaum and Moishezon.
\end{remark}

\bigskip
\noindent \textit{Acknowledgments.} 
Many thanks to Bob Gompf and Danny Ruberman for helpful comments. The author is supported by the NSF Grant DMS-$1510395$.

\end{document}